\documentclass[12pt]{article}
\usepackage[top=2cm, left=3.5cm,right=2cm,bottom=1cm]{geometry}
\usepackage[utf8]{inputenc} 
\usepackage[T1]{fontenc}    
\usepackage{hyperref}       
\usepackage{url}            
\usepackage{booktabs}       
\usepackage{amsfonts}       
\usepackage{nicefrac}       
\usepackage{microtype}      
\usepackage{lipsum}
\usepackage{fancyhdr}       
\usepackage{graphicx}       
\graphicspath{{media/}}     

\usepackage[pages=all, color=black, position={current page.south}, placement=bottom, scale=1, opacity=1, vshift=5mm]{background}
\SetBgContents{
	\tt} 
\usepackage{textcomp}

\usepackage{amsmath, amssymb, amscd, amsthm, amsfonts,eucal,latexsym,dsfont,stmaryrd,enumerate}
\usepackage{amsxtra}
\usepackage{latexsym}
\usepackage{dsfont}
\usepackage{mathrsfs}
\usepackage{niceframe}
\usepackage{color}
\usepackage[all]{xy}
\usepackage{hyperref}
\usepackage{graphicx}
\usepackage{mathtools}
\usepackage{caption}
\usepackage{tikz} 
\usepackage{tikz-cd}  

\usepackage[utf8]{inputenc}
\usepackage{hyperref}
\hypersetup{
	unicode,
	colorlinks,
	breaklinks,
	urlcolor=cyan, 
	linkcolor=blue, 
}

\newcommand{\Spec}{\text{Spec }}

\theoremstyle{plain}
\newtheorem{theorem}{Theorem}[section]

\newtheorem{lemma}[theorem]{Lemma}

\newtheorem{proposition}[theorem]{Proposition}

\theoremstyle{definition}

\newtheorem{remark}{Remark}[theorem]

\usepackage{graphicx, color}
\graphicspath{{fig/}}

\usepackage{algorithm, algpseudocode} 
\usepackage{mathrsfs} 
\usepackage{mathtools}

\renewcommand{\H}{\mathrm{H}}
\newcommand{\h}{\mathrm{h}}

\usepackage{fancyhdr} 
\title{\bf On Some Variants Of Schinzel's Theorem For Global Function Fields
\date{}
}
\author{
  {\normalfont  Nguyen Quang Khai}\\
  {\it Dedicated to the memory of Andrzej Schinzel}
}
\begin{document}
\maketitle
\begin{abstract}
In this paper, we obtain global function field versions of the results of Schinzel and Postnikova for one-dimensional tori, and of  Hahn and Cheon for elliptic curves, which is an analog of the former result.
\end{abstract}
\renewcommand{\thefootnote}{}
\renewcommand{\thefootnote}{\arabic{footnote}}
\setcounter{footnote}{0}
\section{Introduction}
 In 1974, 
Schinzel in \cite{Sch74} proved that for any number field $K$, there exists a constant $n(K)$ depending only on $K$ such that for every $x\in K^*$ which is not a root of unity and every integer $n>n(K)$, the number $x^n-1$ has a primitive divisor.  This is more general than the main result of \cite{Post}. 
Passing to the elliptic curves, S. Hahn and J. Cheon proved a similar result which says that for any elliptic curve $E$ over a number field $K$, and $P$ is a non-torsion $K$-point on $E$, then for every sufficiently large integer $n$, there exists a prime $\mathfrak p$ of good reduction so that the order of $P$ modulo $\mathfrak p$ is equal to $n$. In this paper, we prove the global function field versions for those above results, see Section \ref{the-main-theorem} for elliptic curves and Section \ref{one-dimensional-torus} for one-dimensional tori. We follow a similar strategy in those results, but we need to modify it so that the proof works well in positive characteristics. We note that the result clearly does not hold for unipotent groups.
\section{Some Preliminaries}\label{defn}
 First, we will give proof for the elliptic curve case. In this case, we always assume that the characteristic of the base field is different from 2 and 3.   Let $K$ be a global function field over $\mathbb F_q$ whose constant field is $\mathbb F_q$, where $q=p^s$ for some prime number $p\neq 2,3$.  We denote by  $M_K$ the set of all pairwise non-equivalent discrete
valuations $v$ on $K$ defined as in \cite{LangFDG} Chapter 2, \textsection 3. So for each $v\in M_K$, we have  $\mathfrak p_v$ the prime divisor associated to $v\in M_K$, $\mathcal O_v$ the associated discrete valuation ring with the valuation $v(\cdot)$, and the normalized norm $|\cdot|_v$. By \emph{$\log$} we always mean the logarithmic function with base $q$. We remark that the field $K$ with \textit{multiplicities}  1 satisfies the product formula as in \cite{Zimmer}. So we can use results in \cite{Zimmer}.  For a point $P\in \mathbb P^n(K)$ with homogeneous coordinates $P=[x_0:...:x_n]$,  the \textit{multiplicative height of $P$ relative to $K$}
is $$\H_K(P):=\prod_{v\in M_K}\max\{|x_0|_v,...,|x_n|_v\}. $$
and the \textit{logarithmic height relative to K} is defined as
$$\h_K(P):=\log(\H_K(P))=-\sum_{v\in M_K}\min\{v(x_0),...,v(x_n)\}. $$
Since $p\neq 2,3$, we may assume that an elliptic curve $E$ over $K$  is given by $$y^2=x^3+ax+b\text{ } (a,b\in K)$$ with the identity element $\mathcal O$. For $P=(x,y)\in E(K)$, the \textit{local height function} $\h_v$  at $P$ is  $$ \h_v(P):=\left\{
	\begin{array}{ll}
		-\min\{0,v(x),v(y)\}  & \mbox{if } P \neq \mathcal O \\
		 0& \mbox{if } P=\mathcal O
	\end{array}
\right.$$
The \textit{Weil height} $\h$ at $P$ is defined by $$\h(P):=\sum_{v\in M_K}\h_v(P).$$
Further, for $f\in K(E)$, the \textit{height on $E$ (relative to $f$)} is the function $$\h_f:E(K)\to \mathbb R,\text{ }\h_f(P)=\h_K(f(P)).$$ 
As in \cite{Zimmer},  we have the \textit{N\'eron-Tate} height $\hat{\h}$ on $E(K)$. From \cite{Zimmer} \textsection2 we have that the function $\hat{\h}$ is a non-negative quadratic form on $E(K)$ and the difference $\hat{\h}-\h$ is bounded on $E(K)$. Further, $\hat{\h}(P)=0$ if and only if $P$  is a torsion point.
\begin{proposition}[Northcott-type finiteness theorem]\label{3} With  above notation, for any positive number $B$, the set $$\{P\in E(K):\h(P)<B\}$$  is finite. 
\end{proposition}
\begin{proof}
It follows directly from Northcott's theorem. It also can be proven as follows. Since $E$ is an abelian variety over $K$, $E$ admits a $K/k$-trace $(A,\tau)$ (see Theorem 8 of Chapter VIII,  \textsection 3 in \cite{LangAV}). It means that  $A$ is an abelian variety defined over $k$ and $\tau$ is a $K$-homomorphism 
$$A\times_kK\to E $$  satisfying the usual universal property in the set of pairs of this form.   Applying Theorem 5.3 of Chapter 6, \textsection 5 in  \cite{LangFDG}, the set $\{P\in E(K):\h(P)<B\}$  lies in a finite number of cosets of $A(k).$ Since $A$ is of finite type over $k$ and $k$ is finite, the set $A(k)$ is finite. Hence $\{P\in E(K):\h(P)<B\}$ is also finite.
 \end{proof}
\section{Elliptic curve analogue}\label{the-main-theorem}
Now we are ready to state and prove the main theorem. 
\begin{theorem}\label{1}
    Let $E$ be an elliptic curve over some global function field $K$ of characteristic $p\neq2,3$ and let
$P \in E(K)$ be a non-torsion point. Then for every sufficiently large integer $n$ prime to $p$, there exists a prime $\mathfrak p$ of good reduction so that the order of $P$ in the group of points of $E$ modulo $\mathfrak p$ is equal to $n$. 
\end{theorem}

\begin{remark}
    We note that the second claim of the main theorem in \cite{Cheon-Hahn} does not hold in this setting because Roth's theorem does not hold in positive characteristics.
\end{remark}

The strategy of the proof is similar to the one in \emph{loc. cit.}  We may assume that $E$ is given by the equation $y^2=x^3+ax+b$.  Let $S$ be a  finite set containing all the places at which $E$ has bad reduction, all places dividing $\infty$, and all places at which either $a$ or $b$ has the nonzero valuation, i.e.,
		 $$S=\{v\in M_K:S\text{ has bad reduction at } v\}\cup\{v\in M_K: v(a)\neq0\}\cup\{v\in M_K:v(b)\neq0\}.$$
The set $S$ is finite and we denote by $\#S$ its cardinality. First, we need a result that helps us to detect whether a non-torsion point is trivial or not after taking reduction. 
\begin{lemma}\label{5}
Let $v\in M_K\setminus S$, and let $P$ be a non-torsion point of $E(K)$. Then 
\begin{itemize}
    \item  If $P$ modulo $\mathfrak p_v$ does not equal $\mathcal{O}$, we have $\h_v(P)=0$.
    \item If $P$ modulo $\mathfrak p_v$ equals $\mathcal{O}$, we have $\h_v(nP)=\h_v(P)>0$  for any positive interger $n$ prime to $p$. 
\end{itemize} 
\end{lemma}
\begin{proof}
We may write $P=(x,y)$ and $P=[X:Y:Z]$ in the corresponding projective closure of $E$, where $X,Y,Z\in \mathcal O_K$.   
The condition $P$ mod $\mathfrak p_v=\mathcal{O}$ means that $v(X)>v(Y)$, $v(Z)>v(Y)$, and hence $v(y)<0.$ Therefore, the condition $P$ mod $\mathfrak p_v\neq\mathcal{O}$ is equivalent to either $v(X)\leq v(Y)$ or $v(Z)\leq v(Y)$. If $v(Z)\leq v(Y)$, then $v(y)\geq0$, and from $y^2=x^3+ax+b$ we obtain $v(x)\geq 0$ (since if $v(x)<0$, then $2v(y)=v(x^3+ax+b)=3v(x)<0$, a contradiction), i.e., $\h_v(P)=0$. If $v(X)\leq v(Y)$ and $v(Y)<v(Z)$, then $v(X)<v(Z)$.  But then from the homogeneous Weierstrass equation $Y^2Z=X^3+aXZ^2+bZ^3$ we obtain $$2v(Y)+v(Z)=3 v(X),$$a contradiction. Therefore the first statement is proven. For the second one, from the Weierstrass equation $y^2=x^3+ax+b$, we have $ 3v(x)=2v(y)<0,$ and then $\h_v(P)=-v(y)>0.$ Moreover, if we let $$E_1(K_v):=\{M\in E(K_v):M\text{ mod }\mathfrak p_v=\mathcal{O}\}, $$
we then have an isomorphism of groups (see \cite{Silverman}, Proposition VII.2.2),
$$E_1(K_v)\to \hat{E}(\mathfrak p_v) ,\text{ }M=(x(M),y(M))\mapsto z(M)=\dfrac{-x(M)}{y(M)}, $$ where $\hat{E}$ is the formal group associated to $E.$ Further, this isomorphism gives us the formula $v(y(M))=-3v(z(M)).$ Thus, via this isomorphism, $nP$ maps to $[n].\Big(\dfrac{-x}{y}\Big)=n.\Big(\dfrac{-x}{y}\Big)+\text{(higher-order terms)}$, here $[n].\Big(\dfrac{-x}{y}\Big)$ is $\Big(\dfrac{-x}{y}\Big)+\Big(\dfrac{-x}{y}\Big)+...+\Big(\dfrac{-x}{y}\Big)$ ($n$ times) in $\hat{E}(\mathfrak p_v)$. Since $v(x)>v(y)$, $v(n)=0$ and $\hat{E}(\mathfrak p_v)$ is defined over $\mathfrak p_v$, we obtain $$ v(z(nP))=v\Big([n].\Big(\dfrac{-x}{y}\Big)\Big)=v\Big(\dfrac{-x}{y}\Big)=v(x)-v(y).$$ Consequently, we get $$v(y(nP))=-3v(z([n].P))=3v(y)-3v(x)=3v(y)-2v(y)=v(y)<0.$$ So $3v(x(nP))=2v(y(nP))<0$, and we obtain $$\h_v(nP)=-v(y(nP))\text{, and hence }\h_v(nP)=\h_v(P)>0. $$
\end{proof}
Now, to prove the main theorem, we will give estimates for places in $S$ (as in Lemma \ref{Siegel-cor}) and places in $M_K\setminus S$ (as in the above lemma) and combine them together to deduce a contradition.
\begin{proof}[Prove of Theorem \ref{1}]

	Assume that for any sufficiently large $n>1$ not divisible by $p$, the order of $P$ modulo $\mathfrak p_v$ is not equal to $n$ for any $v\in M_{K}$. In other words, if $nP$ modulo $\mathfrak p_v$ equals $\mathcal O$, then there exists some prime divisor $r$ of $n$ such that $\dfrac{n}{r}P$ modulo $\mathfrak p$ equals $\mathcal O$, and hence, Lemma \ref{5} gives us $$\h_v(nP)=\h_v\Big(\dfrac{n}{r}P\Big).$$
    It follows that for $v\in M_{K}\setminus S$, we have \begin{equation}\label{01}
        \h_v(nP)\leq \sum_{r}\h_v\Big(\dfrac{n}{r}P\Big) ,
    \end{equation}
    where $r$ runs over the set of prime divisors of $n$. For $v\in S$, Lemma \ref{Siegel-cor} yields $$\lim_{n\to\infty}\dfrac{\h_v(nP)}{\h(nP)}=0.$$
    Since $\#S$ is finite, it follows that for any $\epsilon>0$,  
    \begin{equation}\label{02}
        \h_v(nP)\leq\epsilon \h(nP) 
    \end{equation}
    for all sufficiently large  integers $n$. Combining inequalities (\ref{01}) and (\ref{02}), we get
    $$\h(nP)=\sum_{v}\h_v(nP) \leq\sum_{v\not\in S}\sum_{r|n}\h_v\Big(\dfrac{n}{r}P\Big)+\sum_{v\in S}\epsilon. \h(nP)\leq \sum_{r|n}\h\Big(\dfrac{n}{r}P\Big)+\#S.\epsilon. \h(nP).$$
    So \begin{equation}\label{eq2}
        (1-\#S.\epsilon)\h(nP)\leq \sum_{r|n}\h\Big(\dfrac{n}{r}P\Big). 
    \end{equation} 
    Now because there exists a constant $c$ such that $$\hat{\h}(Q)-c<\h(Q)<\hat{\h}(Q)+c\text{ for all }Q\in E(K),$$
    combining this with (\ref{eq2}) implies
 $$(1-\#S.\epsilon)(\hat{\h}(nP)-c) <\sum_{r|n}\hat{\h}\Big(\dfrac{n}{r}P\Big)+c.n $$
 since $\#\{\text{prime divisors of }n\}< n$. Because of the quadraticity of $\hat{\h}$, it follows that 
 $$ (1-\#S.\epsilon)(n^2.\hat{\h}(P)-c)<\sum_{r|n}\dfrac{n^2}{r^2}\hat{\h}(P)+c.n<\dfrac{n^2}{2}\hat{\h}(P)+cn$$ since $\displaystyle\sum_{r|n}\dfrac{1}{r^2}<\dfrac{1}{2}$. Therefore $$\Big(\dfrac{1}{2}-\#S.\epsilon\Big)n^2.\hat{\h}(P)<(n+1-\#S.\epsilon).c $$ We choose $\epsilon<\dfrac{1}{2\#S}$ and let $n$ tend to $\infty$, then we obtain $\hat{\h}(P)=0$, a contradiction.
\end{proof}
\begin{remark}

We also need the following result that is similar to a classical theorem of Siegel for number fields, see \cite{Silverman} Theorem 3.1.
\begin{lemma}\label{Siegel-cor}
Let $E/K$ be an elliptic curve over a global function field of characteristic $p\neq 2$ and 3, and let $P\in E(K)$ be a non-torsion point.  Then for any place $v\in M_K$, we  have
 $h_v(nP)$ is bounded above when $n$ ranges over integers coprime to $p$. In particular, we have
$$\lim\limits _{\substack{%
  \gcd(n,p)=1  \\n\to\infty
     }}\dfrac{h_v(nP)}{h(nP)}=0.$$
\end{lemma}
\begin{proof}
Since the goal is to obtain the boundedness, replacing $x$ by $u^2x$ for some $u\in K^\times$, we may assume that the given Weierstrass equation is minimal at $v$. The desired boundedness then follows from the following observation: if there are $n\neq m$ and $n,m$ coprime to $p$ such that both $nP$ and $mP$ mod $\mathfrak p_v$ equal $\mathcal{O}$, then 
$h_v(nP)>0$, $h_v(mP)>0$ and $h_v(nP)=h_v(nmP)=h_v(mP)$ by Lemma \ref{5}. 


Since $h(nP)$ tends to $\infty$ as $n$ tends to $\infty$ (thanks to the Northcott property), the second claim follows.
\end{proof}

We note that the condition $\mathrm{gcd}(n,p)=1$ is necessary. For example, consider a supersingular elliptic curve $E$ over $K$. Equivalently, the Hasse invariant $A(E,\omega)$ is equal to 0, where $\omega=\dfrac{dx}{y}$, see \cite{Katz-Mazur}. We note that when $p>3$, $A(E,\omega)$ is the coefficient of $x^{p-1}$ in $(x^3+ax+b)^{(p-1)/2}.$ Therefore, for $v\not\in S$ where $S$ is the set in the previous proof, the reduction modulo $\mathfrak p_v$ of $E$ is an elliptic curve $E_{v}$ whose Hasse invariant is also 0. Thus  $E_{v}$ is supersingular. So for such $v$, $E_v[n]=E_v[np]$ for all integer $n$. Therefore, for all but finitely many $v$, the order of the reduction of the point $P$ modulo $\mathfrak p_v$ must be prime to $p$.
\end{remark}
When $E$ is ordinary, we have 
\begin{theorem}\label{ordinary-ec}
    Let $E$ be an ordinary elliptic curve over some global function field $K$ of characteristic $p\neq2,3$ and let $P \in E(K)$ be a non-torsion point. We fix a positive integer $t$. Then for every sufficiently large integer $n$ prime to $p$, there exists a prime $\mathfrak p$ of good reduction so that the order of $P$ in the group of points of $E$ modulo $\mathfrak p$ is equal to $np^t$.  
\end{theorem}
\begin{proof}
Because $E$ is ordinary, there exists a point $Q\in E(\Bar{K})$ of order $p^t$. We set $L:=K(E[p^t])$ the $p^t$-division field of $E$, and denote $E_L$ the base change of $E$ to $L$. Then $P-Q\in E_L(L)$ is also a point of infinite order. We note the following properties of the reduction of points (for more details on reductions and integral models, we refer to \cite{Perucca}). 
\begin{enumerate}
    \item (see \cite{Perucca} Lemma 1.2.3) For all but finitely many primes $\mathfrak p$ of $K$ the following holds:  for any prime $\mathfrak q$ of $L$ lying over $\mathfrak p$, the order of $P$ (as an $L$-point of $E_L$) modulo $\mathfrak q$  equals the order of $P$ modulo $\mathfrak p$.
    \item (see \cite{Perucca} Corollary 2.3.4) Since $Q$ is torsion, the order of $Q$ modulo $\mathfrak q$ equals the order of $Q$, which is $p^t$, for all but finitely many primes $\mathfrak q$ of $L$. 
\end{enumerate}
We note that the proofs of those properties given
in \cite{Perucca} also work well over global function fields. We call $V$ the set of  {\it exceptional primes} of $K$ in (1)
and  call $U$ the set of {\it exceptional primes} of $L$ in (2)
and primes of $L$ lying above primes in $V$. Then both $V$ and $U$ are finite. Now we apply Theorem \ref{1} for $P-Q\in E_L(L)$, we have for every sufficiently large integer $n$ prime to $p$, there exists a prime $\mathfrak q$ of good reduction so that the order of $P-Q$ modulo $\mathfrak q$ equals $n.$ Since $U$ is finite, the prime $\mathfrak q$  does not lie in $U$ for $n$ sufficiently large. Since $np^tP=np^t(P-Q)+np^tQ=np^t(P-Q)$,  the point $np^tP$ modulo $\mathfrak q$ equals $\mathcal{O}$. So the order of $P$ is of the form $mp^s$ where $m|n$ and $s\leq t$. Then $$\mathcal{O}=mp^tP\mod\mathfrak q=mp^tQ+mp^t(P-Q)\mod\mathfrak q=mp^t(P-Q)\mod\mathfrak q, $$
and hence, $n|mp^t$ which implies that $m=n$. Similarly, we have $$\mathcal{O}=np^sP\mod\mathfrak q=np^sQ+np^s(P-Q)\mod\mathfrak q=np^sQ\mod\mathfrak q. $$Thus $p^t|np^s$ which means that $s=t$. Therefore the order of $P$ modulo $\mathfrak q$ is $np^t$. Since  $\mathfrak p$,  the prime of $K$ lying under  $\mathfrak q$,  does not lie in $V$, the order of $P$ modulo $\mathfrak p$  is also equal to $np^t$. The theorem is then proven.
 \end{proof}
\section{One-dimensional torus version}\label{one-dimensional-torus}
In this section, we prove an analogous result for one-dimensional tori. First, we establish an analogous result for the multiplicative group $\mathbb G_m$ over $K$ in  arbitrary characteristic (not necessarily different from 2 and 3).
\begin{proposition}\label{multiplicative-group}
    Let $K$ be a global function field, and let $x\in K\setminus\{0\}$ be not a root of unity. Then for every positive integer $n>1$ prime to $p$, there exist a place $v\in M_K$ such that the order of the reduction of $x$ modulo $\mathfrak p_v$ in the group $\mathbb G_m(\mathcal{O}_v/\mathfrak p_v)$ is equal to $n$.   
\end{proposition}
\begin{remark}
The condition $n$ being prime to $p$ is necessary. Indeed, for any $x\in K$, we have $x^{np}-1=(x^n-1)^p$. Therefore, if $\mathfrak p_v|x^{np}-1$ for some place $v$, then $\mathfrak p_v|x^n-1$.
\end{remark}
Now we prove Proposition \ref{multiplicative-group}
\begin{proof}[Proof of Proposition \ref{multiplicative-group}]
We denote by $\mu$ the M\"obius function.  Let  $W$ be the set $\{v\in M_K:v(x^n-1)>0\}.$ We consider four following cases. 
\begin{enumerate}
    \item $v\in W$ such that $n$ is not the order of $x$ modulo $\mathfrak p_v$, we call this order by $n_0.$   Then $n=n_0.k$ for some positive integer $k$ and
 $$x^n-1=(x^{n_0}-1)(x^{n_0(k-1)}+x^{n_0(k-2)}+...+x^{n_0}+1). $$
 Since $x^{n_0(k-1)}+x^{n_0(k-2)}+...+x^{n_0}+1\equiv k\not\equiv0\mod \mathfrak p_v$ ($(k,p)=1$), we have $v(x^n-1)=v(x^{n_0}-1)$. Thus 

\begin{equation} 
\begin{split}
v(\Phi_n(x)) & =\sum_{m|n}\mu\Big(\dfrac{n}{m}\Big)v(x^m-1)   = \sum_{n_0|m|n}\mu\Big(\dfrac{n}{m}\Big)v(x^m-1)\\
 & = \sum_{n_0|m|n}\mu\Big(\dfrac{n}{m}\Big)v(x^{n_0}-1)=0\text{ since }n_0<n.
\end{split}
\end{equation}
\item $v\in M_K$ satisfying $v(x)>0$. Then $v(x^m-1)=0$ for all positive integer $m$. It implies that $v(\Phi_n(x))=0.$
\item $v\in M_K$ satisfying $v(x)<0$. Then $v(x^m-1)=v(1-x^{-m})+v(x^m)=mv(x).$ Hence 
\begin{align*}
v(\Phi_n(x))  =\sum_{m|n}\mu\Big(\dfrac{n}{m}\Big)v(x^m-1)   
  = \sum_{m|n}\mu\Big(\dfrac{n}{m}\Big)mv(x) = \phi(n).v(x).
\end{align*}
\item $v\in M_K$ satisfying $v(x)=0$ and $v\not\in P$. Then $v(x^m-1)=0$ for all $m|n$, and hence $v(\Phi_n(x))=0$.
\end{enumerate}
Combining these equalities, we obtain that if for every $v\in W$, $n$ is not the order of $x$ modulo $\mathfrak p_v$, then we have 
\begin{align*}
0=\sum_{v\in M_K}v(\Phi_n(x))  = \sum_{v\in M_K:v(x)<0}v(\Phi_n(x)) = \sum_{v\in M_K:v(x)<0}\phi(n).v(x).
\end{align*}
This equality holds if and only if there is no place $v$ such that $v(x)<0$, which means that $x$ must lie in the constant field $\mathbb F_q$, which is a contradiction. Thus, there exists some $v$ in $P$ such that $n$ is the order of $x$ modulo $\mathfrak p_v$, this is what we want.
\end{proof}
To prove the same result for an arbitrary one-dimensional torus, we just extend it to a multiplicative group and use some arguments on reductions.
\begin{theorem}\label{tori}
    Let $G$ be a one-dimensional torus over a global function field $K$, and let $x\in G(K)$ be a point of infinite  order. Then for every sufficiently large integer $n$ prime to $p$, there exists a place $v\in M_K$ such that the reduction of $x$ modulo $\mathfrak p_v$ exists, and the order of $x$ modulo $\mathfrak p_v$ is $n.$  
\end{theorem}

\begin{proof}
This theorem can be proven using the method of Lemma 1.2.3 in \cite{Perucca}.  Firstly, after discarding a finite set of places $S$, we may assume that the coefficients of $x$ are in $\mathcal{O}_S$, the ring of $S$-integers in $K.$ After discarding finitely more  places (we still denote this set of places by $S$),   we can assume that $G$ admits  an integral model  $\underline{G}$ over $\mathcal{O}_S$ and $x\in\underline{G}(\mathcal{O}_S)$. Further, there exists a finite Galois extension $K'$ of $K$ such that $G\times_KK'=\mathbb G_{m,K'}$. After discarding finitely more places (we still denote this set of places by $S$), one can assume that $\mathbb G_{m,K'}$ admits an integral model $\mathbb G_{m,\mathcal{O}_S'}$ ($S'$ is the set of places in $M_{K'}$ above $S$) which is an extension $\underline{G}$, i.e., $\mathbb G_{m,\mathcal{O}_{S'}}=\underline{G}\times_{\mathcal{O}_S}\mathcal{O}_{S'}$ (thanks to the uniqueness of integral models). Now for any $v\not \in S$ and $w\not\in S'$ above $v$, let $P$ be an $\mathcal{O}_v$-point of $\underline{G}$.  Then $P$ is also an $\mathcal{O}_v$-point of $G$, and it can be lifted to a point $\mathcal{O}_w$-point $P'$ of $\mathbb G_{m,K'}$, since $\mathcal{O}_v\otimes_KK'\cong\displaystyle\prod_{w|v}\mathcal{O}_w$. Since $\mathcal{O}_{S'}\subset \mathcal{O}_w$, the point $P$ is also an $\mathcal{O}_{w}$-point of $\mathbb G_{m,S'}$, i.e., we have the commutative diagram 
\[\begin{tikzcd}
	{\Spec\mathcal O_w} & {\mathbb G_{m,\mathcal O_{S'}}} & {\mathbb G_{m,K'}} \\
	{\Spec \mathcal O_v} & {\underline G} & {G}
	\arrow["{P'}", from=1-1, to=1-2]
	\arrow[from=1-2, to=1-3]
	\arrow[from=1-1, to=2-1]
	\arrow["{P}", from=2-1, to=2-2]
	\arrow[from=2-2, to=2-3]
	\arrow[from=1-2, to=2-2]
	\arrow[from=1-3, to=2-3]
\end{tikzcd}\]
In other words, we have $\underline G(\mathcal O_v)\subset \mathbb G_{m,\mathcal O_{S'}}(\mathcal O_w)$. Taking reduction, we see that the group of reduction points modulo $\mathfrak p_v$ in $\underline G$ is injected in the  group of reduction points modulo $\mathfrak p_w$ in $\mathbb G_{m,\mathcal O_{S'}}$, and hence the order of $P$ modulo $\mathfrak p_v$ is equal to the order of $P'$ modulo $\mathfrak p_w$ for every $w|v$ and every $P\in \underline{G}(\mathcal{O}_v)$. Now we take $P$ to be $x\in \underline G(\mathcal{O}_S)$, then $x$ lifts to $x'\in\mathbb G_{m,\mathcal{O}_{S'}}(\mathcal{O}_{S'})=\mathcal{O}_{S'}^\times\subset K$, and we have the following commutative diagram
\[\begin{tikzcd}
	{\Spec\mathcal O_w} & {\Spec\mathcal O_{S'}} & {\mathbb G_{m,\mathcal O_{S'}}} \\
	{\Spec \mathcal O_v} & {\Spec\mathcal O_{S}} & {\underline G}
	\arrow[from=1-1, to=1-2]
	\arrow["x'", from=1-2, to=1-3]
	\arrow[from=1-1, to=2-1]
	\arrow[from=2-1, to=2-2]
	\arrow["x", from=2-2, to=2-3]
	\arrow[from=1-2, to=2-2]
	\arrow[from=1-3, to=2-3]
\end{tikzcd}\]
 Applying Proposition \ref{multiplicative-group}, for every sufficiently large integer $n$ prime to $p$, since $S'$ is finite, there exists  a place $w_0\not\in S' $ such that the order of the reduction of $x'$ modulo $\mathfrak p_{w_0}$ equals $n$. Thus the order of the reduction of $x$ modulo $\mathfrak p_{v_0}$, with $v_0\not\in S$ is the place that lies under $w_0$,  also equals $n$.
\end{proof}

\section*{Acknowledgments}
I would like to express my deepest gratitude to Professor Nguyen Quoc Thang for his supervision of this project and his support. I would like to thank Professor Yann Bugeaud for indicating that Roth's theorem does not hold in positive characteristics. This work is funded by the International Center for Research and Postgraduate Training in Mathematics under the auspices of UNESCO, grant ICRTM03$\_$2020.06, and supported by the Domestic Master Scholarship Programme of Vingroup Innovation Foundation, Vingroup Big Data Institute, grant VINIF.2021.ThS.05.


{Institute of Mathematics, Vietnam Academy of Science and Technology, Hanoi, Vietnam}\\
\indent Email: \texttt{nqkhai@math.ac.vn} 


\end{document}